\documentclass[journal]{IEEEtran}
\usepackage{amsmath,amsfonts}
\usepackage{algorithmicx}
\usepackage{algorithm}
\usepackage{array}
\usepackage[caption=false,font=normalsize,labelfont=sf,textfont=sf]{subfig}
\usepackage{textcomp}
\usepackage{stfloats}
\usepackage{url}
\usepackage{verbatim}
\usepackage{graphicx}
\usepackage{cite}
\usepackage{supertabular,booktabs}
\usepackage{comment}
\usepackage{xcolor}
\usepackage{longtable}
\usepackage{eurosym}
\usepackage{comment}
\usepackage{algpseudocode}
\begin{document}

%\title{Multi-Objective Optimization for Energy Communities: Balancing Cost
%and Resilience Considering the Impact of Energy Poor Members}
\title{Cost versus Resilience in Energy Communities: A Multi-Objective Member-Focused Analysis}
\author{Lia Gruber,~Sonja Wogrin
        % <-this % stops a space
% <-this % stops a space
\thanks{
Working paper submitted March 2025 (Corresponding Author: Lia Gruber) \\
L. Gruber and S. Wogrin are with the Institute of Electricity Economics and Energy Innovation, Graz
University of Technology, Graz, Austria.\\
\\(e-mail: lia.gruber@tugraz.at; ; wogrin@tugraz.at)}}

\markboth{Working Paper}
{L. Gruber, S. Wogrin \MakeLowercase{Cost versus Resilience in Energy Communities: A Multi-Objective Member-Focused Analysis}}

\maketitle

\begin{abstract}
This paper develops a multi-objective optimization framework to analyze the trade-offs between annual costs and resilience in energy communities. Under this framework, three energy community operation strategies are analyzed: a reference case where all assets are member-owned, implementing a communal battery electric storage system, and subsidizing energy-poor members. The results indicate that increasing resilience leads to higher operational costs and smaller feasible ranges of energy community energy prices. The analysis reveals that those trade-offs have a heterogeneous impact across different member groups. Owners photovoltaics are most affected due to curtailed energy. Notably, the study shows that while implementing community-owned storage does not directly provide financial benefits to energy-poor members, alleviating the energy price for these members leads to an overall cost reduction of more than 30\%. This research provides insights into the operational complexity of energy communities and highlights the importance of technologically robust and socially inclusive energy communities.

\end{abstract}
%\begin{comment}
\section*{Nomenclature}
\label{sec:Nomenclature}
\subsection{Indices}
\begin{IEEEdescription}[\IEEEusemathlabelsep\IEEEsetlabelwidth{$a(m,n,l)$}]
\item[$p$] Type of objective function
\item[$h$] Time periods in hours
\item[$n$] Members
\item[$g$] Generating units
\item[$s(g)$] Subset of storage generation units
\item[$r(g)$] Subset of renewable generation units
\item[$e(g)$] Subset of EVs
\item[$i$] Bus in electricity grid
\item[$itrans(i)$] Transmission node 
\item[$gc(g)$] Subset of community owned generation units
\item[$gn(g,n)$] Subset of generation units owned by members
\end{IEEEdescription}

\subsection{Parameter}
\label{subsec:Parameter}
\begin{IEEEdescription}[\IEEEusemathlabelsep\IEEEsetlabelwidth{$a(m,n,l)$}]
\item[$RES^{max}$] Maximum of resilience (MW)
\item[$RES^{min}$] Minimum of resilience (MW)
\item[$gp$] Grid point 
\item[$C^{gt,flat}$] Grid tariff fixed costs (\euro/member)
\item[$C^{ovh,acc}$] EC overhead cost accounting (\euro)
\item[$C^{ovh,met}$] EC overhead cost per metering point (\euro/metering point)
\item[$C^{ovh,var}$] EC overhead cost per consumption (\euro/MWh)
\item[$C_{h}^{import}$] Cost of import (\euro/MWh)
\item[$C_{h}^{export}$] Cost of export (\euro/MWh)
\item[$C^{gt,flat}$] Grid tariff fixed costs (\euro/member) 
\item[$C^{gt,in}$] Grid tariffs for EC energy (\euro/MWh) 
\item[$C^{gt,out}$] Grid tariffs for imported energy (\euro/MWh) 
\item[$Cost_n^{in}$] Member cost inside the EC (\euro)
\item[$Cost_n^{out}$] Member cost outside the EC (\euro)
\item[$Cost^{c,in}$] EC cost inside the EC (\euro)
\item[$Cost^{c,out}$] EC cost outside the EC (\euro)
\item[$Cost^{c,ovh}$] EC overhead cost (\euro)
\item[$Cost^{c,net}$] EC cost balance (\euro)
\item[$Cost^{c,net,ep}$] EC cost balance with energy poverty subsidy(\euro)
\item[$C^{EC,g}$] EC energy price for producers (\euro/MWh)
\item[$C^{EC,n}$] EC energy price for consumers (\euro/MWh)
\item[$D_{h,n}$] Member demand (MW)
\item[$D^c_{h}$] EC demand (MW)
\item[$D_{h,n}^{net}$] Member net demand (MW)
\item[$P_{h}^{exe,n}$] Member excess energy (MW)
\item[$P_{h}^{c,net}$] EC net energy production (MW)
\item[$P_{h}^{c}$] EC energy production (MW)
\item[$Rd^c_h$] EC residual grid demand (MW) 
\item[$Rd_{h,n}$] Member residual grid demand (MW) 
\item[$Exe_h^c$] EC Excess energy (MW)
\item[$Exe_{h,n}^{out}$] Member excess energy sold in EC (MW)
\item[$Exe_{h,n}^{in}$] Member excess energy sold outside EC (MW)
\item[$Exe_h^{gc}$] Excess energy from communal assets (MW)
\item[$P_h^{dis,c}$] EC actual distributed energy (MW)
\item[$Share_{h,n}^{d,\%}$] Dynamic share of production (\%)
\item[$Sc_{h,n}$] Member Self-consumption (MW)
\end{IEEEdescription}  
\subsection{Variables}
\label{subsec:Variables}
\begin{IEEEdescription}[\IEEEusemathlabelsep\IEEEsetlabelwidth{$a(m,n,l)$}]
\item[$s_p$] Slack variable
\item[$r_p$] Range of objective function
\item[$\epsilon_p$] Epsilon of objective function
\item[$TC$] Total cost (\euro)
\item[$RES$] Resilience (MW) 
\item[$imp_{h,itrans}$] Import (MW)
\item[$exp_{h,itrans}$] Export (MW)
\item[$cs_{h,g}$] Consumption (MW)
\item[$p_{h,g}$] Production (MW)

\end{IEEEdescription}  
%\end{comment}
\newpage
\section{Introduction}
\subsection{Motivation}
\IEEEPARstart{R}{enewable} Energy Communities (RECs), defined in the European Union's Renewable Energy Directive II (REDII), provide the regulatory basis for collective energy generation, consumption, storage, and selling \cite{REDII}. The directive also gives a direction for overlying goals that RECs should have. Their primary purpose should not be financial gain but the provision of environmental, economic, or social benefits to the energy community (EC) members and/or the local area around the community. This development is part of a paradigm shift in the energy sector, where citizens are empowered to take a more active role. The technological promise of RECs lies in their potential to increase the penetration of renewables in the grid. By creating an incentive to consume energy geographically close to where it is produced, one hope is that RECs could reduce the strain that the increase of intermittent renewables can have on the grid, enable more renewable capacity, and lower overall system costs by decreasing the need for grid reinforcements.

Beyond technical benefits, energy communities can potentially deliver significant social value by aiding energy poverty. The European Commission has been addressing energy poverty since its third energy package (Directive 2009/72/EC \cite{EUDir72} and Directive 2009/73/EC, \cite{EUDir73}), which introduced the concept of energy poverty but left the responsibility for its definition and mitigation to individual Member States. The REDII framework acknowledges the role of RECs in fighting energy poverty by enabling lower energy prices.\cite{Doukas2020} However, depending on how the operation of the EC is organized, there are potential barriers to participation, such as high initial investment costs into assets that prevent vulnerable households from accessing the benefits of ECs. Therefore, mechanisms to support energy-poor members are critical to ensuring that ECs deliver broad-based social benefits and are accessible to all ranges of society.

\subsection{Related Work and Contributions}

\subsubsection{EC impact on the distribution grid}
Beyond cost and environmental objectives, ECs are increasingly viewed as tools for improving grid management and enabling greater penetration of renewables. The ability of ECs to provide grid services through the strategic operation of flexibility assets and demand-side management has been the focus of a lot of research. Polgari et al. in \cite{Polgari2024} demonstrated that control of voltage and reactive power (P(U), Q(U)) within ECs can increase Photovoltaic (PV) capacity without running into voltage limits. In a case study of an EC in Hungary, they show that it is possible to enable more renewable generation without overloading the distribution grid. Sudhoff et al. in \cite{Sudhoff2022} showed that grid-friendly EC operation, driven by economic incentives, can reduce peak load at the low-voltage substation by 23–55\%. Dimovski et al. in \cite{DIMOVSKI2023} used Monte Carlo simulations to assess how different EC configurations impact grid infrastructure, highlighting the importance of balancing load and production schedules for reducing grid losses and voltage violations. Berg et al. in \cite{BERG2023} studied the role of community-owned batteries in providing services to the distribution grid. Their findings show that financial incentives for peak load reduction could improve grid resilience without increasing costs for EC members. These studies demonstrate that ECs can provide valuable grid services, but this can only be achieved with appropriate operational incentives and financial compensations. 

\subsubsection{Energy Poverty in ECs}
The potential of ECs to reduce energy poverty is not only recognized in policy but also in research. Hanke and Guyet in \cite{Hanke2023} analyzed the role of ECs in enhancing energy justice, concluding that enabling frameworks with clear taxonomy are needed to support ECs in this endeavor. Bielig et al. in \cite{Bielig2022} explored the social impact of ECs via a systemic literature review of studies that measure this impact. They conclude there is also a significant need for more evidence collection on this topic. Boeri et al. in \cite{Boeri2020} indicate that ECs can support vulnerable households by helping them to make their energy demand more efficient. Konstantopoulos et al. in \cite{Konstantopoulos2023} listed a range of actions that Greek ECs can take to address energy poverty, including targeted tariff reductions and direct financial support. Fina in \cite{Fina2024} proposed a redistribution mechanism of EC-incurred grid costs based on income, reducing the financial burden on poor energy households. Ceglia et al. in \cite{Ceglia2022} introduced a social energy poverty index that evaluates the mitigation potential of an Italian EC and puts it at 12-16\%. While these studies show the potential of ECs to reduce energy poverty, none have analyzed how adding a focus on this social dimension affects EC operation and pricing strategies.

\subsubsection{Multi-objective optimization of ECs}
The operational complexity of ECs often results in multiple technical, economic, social, or environmental goals that can be followed and they do not always have to align. Multi-objective optimization (MOO) has been widely applied in EC research to address the complex trade-offs between multiple objectives. Petrelli et al. in \cite{PETRELLI2021} applied a multi-objective, multi-year planning method to microgrids, balancing socio-economic benefits (net present cost, job creation), security (public lighting coverage), and environmental impacts (carbon emissions, land use). Schram et al. in \cite{Schram2020} proposed a community energy storage system to reduce costs and emissions, showing that both can be reduced simultaneously. Ahmadi et al. in \cite{Ahmadi2022} focused on minimizing grid dependency and operational costs within a microgrid. Similarly, Fan et al. in \cite{Fan2022} designed a near-zero energy community using a multi-objective collaborative optimization method to minimize annual carbon emissions, total costs, and grid interactions. Fangjie et al. in \cite{Fangjie2023} evaluated the effect of different demand responses in a community virtual cloud power plant. With a multi-objective grey wolf optimization algorithm (a type of swarm intelligence optimization algorithm), they optimize the energy supplier profit, resident cost, renewable energy utilization rate, and carbon treatment amount. Lazzari et al. in \cite{Lazzari_2023} applied a combinatorial optimization framework for participant selection. In a second step, they minimized PV excess while maximizing profitability for all members. However, the impact on different member groups and EC energy price setting has not been fully integrated into existing models.\\

\subsubsection{Contributions}
The aim of this paper is to analyze how the trade-off between cost and grid resilience in ECs with different operation strategies affects the different member groups. The framework also evaluates how different operation strategies themselves impact the member groups. Therefore, we propose a multi-objective optimization framework to minimize annual EC costs and maximize resilience by minimizing the power peaks at the transformer. The Pareto front is established by using the augmented $\epsilon$-constraint (AUGMECON) method \cite{MAVROTAS2009}. A new ex-post energy allocation mechanism considering different asset ownership options is introduced. It is followed by a cost analysis of the EC and its members, using feasible EC energy price ranges resulting from the Pareto front. The framework is used on a test EC under different operation strategies: only member-owned assets, having a communal BESS, and subsidizing poor energy members. 

The remainder of the paper is organized as follows: Section \ref{sec:Meth} presents the methodology, including the details of the optimization framework and the ex-post allocation mechanism. Section \ref{sec:CS} describes the main characteristics of the test REC and the case studies with their different EC operation strategies including the subsidy of energy poor members. Section \ref{sec:Results} discusses the results, focusing on the Pareto frontier, the cost-resilience trade-off, and the impact on different member groups with a special focus on energy poor members. Section \ref{sec:conclusion} concludes this paper.

\section{Methodology}
\label{sec:Meth}
The optimization model used in this paper is the EC version of the Low-carbon Expansion Generation Optimization Model (LEGO \cite{Wogrin_2022}, EC LEGO \cite{GRUBER2024109592}). To optimize the operation of the EC, a mixed integer problem is solved in GAMS using the Gurobi solver. The model has a high degree of temporal flexibility. In this case, one year with chronological hours was chosen. The EC and its members are modeled via their demand, generation, and flexibility. In this paper, PV is used for the generation units of the EC, and the flexibility options are BESS and EVs. The EC is embedded in a distribution grid the constraints of which are modeled with an optimal power flow. The EC can import and export energy via its transmission node from higher grid levels. In this paper, the EC is based on the low-voltage grid with all members connected behind one transformer that acts as the transmission node $i_{trans}$. The energy and cost allocation between members are calculated ex-post. Those allocation processes are modeled after the Austrian regulatory framework.

\subsection{Objective functions}
\label{subsec:OBF}
Different economic and technological objective function options are implemented in the EC LEGO model. The two used in this study are shown in Eqn. (\ref{eqn:OBCost}) and (\ref{eqn:resil}). The first objective function is the community's annual total cost. This includes costs from energy imports, profits from energy exports, grid tariffs, and overhead costs for managing the EC. 
\begin{align}
\label{eqn:OBCost}
TC = C^{gt,flat}\cdot card(n) + C^{ovh,acc} + \nonumber \\ C^{ovh,met}\cdot card(n) \cdot card(g) 
+\sum_{h} \Big(   \nonumber \\
 +(C^{import}_{h}+C^{gt,out})\cdot imp_{h,itrans}-C^{export}_{h}exp_{h,itrans}\nonumber \\
 +C^{gt,in}\cdot (\sum_{n}D_{h,n}+\sum_{s,e\in g}cs_{g}- imp_{h,itrans})\nonumber \\
\nonumber \\
+ C^{ovh,var}\cdot ( \sum_{n} D_{h,n} +\sum_{s,e\in g} cs_{h,g})
%\nonumber \\
%+ \sum_i C^{ENS} pns_{h,i}+\sum_i C^{EEN} ep_{h,i}\Big)
\Big)
\end{align}
The grid tariff structure is based on the Austrian regulatory framework. It consists of several components, some technical, e.g., for grid use and loss, and others tax, like an electricity duty and a green energy surplus charge. For the promotion of ECs, the tariffs are reduced technological components and exemptions for tax components. For simplicity, all components paid per MWh are presented in the formulation as one coefficient $C^{gt,out}$ for the energy from outside the EC and $C^{gt,in}$ for the energy shared inside the EC. For a more detailed description of the component's structure, see \cite{GRUBER2024109592}. The components paid annually, independent of the energy consumed, are represented by $C^{gt,flat}$. A new addition to the model is the overhead costs of the EC. The billing process and member management of ECs are associated with costs. In Austria, this is often outsourced to specialized service providers. The labor of these processes depends on the size and how active the EC is. Therefore, the overhead cost formulation was split between fixed annual accounting costs $C^{ovh,acc}$,  costs dependent on metering points $C^{ovh,met}$, and variable costs $C^{ovh,var}$ dependent on the consumed energy. 

This study's second objective is to maximize resilience by minimizing the 
power peaks at the transmission node. This is defined by the hourly import and export, as illustrated in Eqn. (\ref{eqn:resil}).
%\begin{equation}
%\label{eqn:OBRes}
%    min~RES
%\end{equation}
\begin{equation}
\label{eqn:resil}
    RES \geq imp_{h,itrans}+exp_{h,itrans}
\end{equation}

\subsection{Multi-Objective Optimization}

This paper reformulated the EC LEGO single-objective model as a multi-objective problem. A generic multi-objective optimization problem is shown in Eqn. (\ref{eqn:MMP}). Here $f_1(\textbf{x})$ to $f_p(\textbf{x})$ are the multiple objective functions where $p\geq 2$. 
They are subject to the set of decision variables \textbf{x} in domain $S$. 
\begin{align}
\label{eqn:MMP}
\max_{\textbf{x}}(f_{1}(\textbf{x}),f_{2}(\textbf{x}),...f_{p}(\textbf{x}))\nonumber \\
s.t. \nonumber \\
\textbf{x} \in S
\end{align}
One way to solve this problem is the $\epsilon$-constraint method, where the multi-objective optimization is transformed into several single-objective optimization models. This is done by using the first objective function as the primary objective function and adding the others as constraints. They are limited by $\epsilon$, as presented in Eqn. (\ref{eqn:eps}). The first step of the method is to calculate the range of the objective functions. While the best value is easy to find, the worst value is more challenging to determine. This range is often produced by creating a payoff table, where the minimum of each column approximates the worst value. However, individual optimization models may not consistently deliver Pareto-optimal solutions, mainly if alternative optima exist. To address this problem, \cite{MAVROTAS2009} developed the augmented $\epsilon$-constraint (AUGMECON) method, where lexicographic optimization ensures that only pareto-optimal solutions are included in the payoff table. 
\begin{align}
\label{eqn:eps}
\max_{\textbf{x}}(f_{1}(\textbf{x}))\nonumber \\
s.t. \nonumber \\
f_{2}(\textbf{x}) \geq \epsilon_2 \nonumber \\
... \nonumber \\
f_{p}(\textbf{x}) \geq \epsilon_p \nonumber \\
\textbf{x} \in S
\end{align}
Another problem with the original $\epsilon$-constraint method is that the solutions are only efficient if all objective function constraints are binding. In the case of alternative optima, the solution may be weakly efficient. AUGMECON resolves this issue by converting the inequality constraints of the secondary objective functions into equalities by adding slack variables. These variables are then included in the objective function, ensuring that only efficient solutions are obtained. This is presented in Eqn. (\ref{eqn:epscm}), where $s_p$ are the slack variables corresponding to each secondary objective function, and $\delta$ is a small scalar to give them a lower priority. In order to avoid scaling problems, the slack variables are divided by the range $r_p$ of their corresponding objective function.

\begin{align}
\label{eqn:epscm}
\max_{\textbf{x},s_p}(f_{1}(\textbf{x})+\delta \cdot (s_2/r_2 + s_3/r_3 + ... + s_p/r_p))\nonumber \\
s.t. \nonumber \\
f_2(\textbf{x})-s_2 = \epsilon_2 \nonumber \\
f_3(\textbf{x})-s_3 = \epsilon_3 \nonumber \\
f_p(\textbf{x})-s_p = \epsilon_p \nonumber \\
... \nonumber \\
\textbf{x} \in S~and~s_p \in \mathbb{R}^+
\end{align}

The objective functions of total annual EC cost and resilience, presented in Eqn.s (\ref{eqn:OBCost}) and (\ref{eqn:resil}), are optimized with the described AUGMECON method. Eqn. (\ref{eqn:epscmCR}) shows the corresponding multi-objective optimization problem. The epsilon of the resilience objective function $\epsilon_{res}$ is calculated in Eqn. (\ref{eqn:epsresil}).
$RES^{max}$ is the maximum resilience, and $RES^{min}$ is the minimum resilience resulting from the optimizations done for building the payoff table. The range $r_{res}$ (Eqn. (\ref{eqn:rangeres})) between $RES^{max}$ and $RES^{min}$ is divided into a number of $gp^{max}$ intervals. $gp$ are the grid points along this range.  
\begin{equation}
\label{eqn:rangeres}
 r_{res} = RES^{max}-RES^{min}
\end{equation}
\begin{equation}
\label{eqn:epsresil}
 \epsilon_{res} = RES^{max} - \frac{gp \cdot r_{res}}{gp^{max}}    
\end{equation}
\begin{align}
\label{eqn:epscmCR}
\max_{\textbf{x},s_p}(-TC+\delta \cdot \frac{s_{res}}{r_{res}}) \nonumber \\
s.t. \nonumber \\
 RES+s_{res} = \epsilon_{res} \nonumber \\
\textbf{x} \in S~and~s_p \in \mathbb{R}^+
\end{align}
\subsection{Constraints}

The production, flexibility, and power flow constraints of the EC-LEGO are detailed in \cite{GRUBER2024109592}, but for the sake of completeness, a short description is presented in the following. Renewables are modeled with hourly capacity factor profiles. The model includes curtailment options for renewables. This type of PV production modeling based on measurements is relatively simple compared to detailed analytical models. Another example of EC optimization models using capacity factors for modeling PV production can be found in \cite{ASKELAND2021}. The formulation of BESS includes state-of-charge (SOC) and simultaneous charging and discharging prohibiting constraints. This BESS formulation can be considered standard and can be found modeled similarly in \cite{VANDERSTELT2018},\cite{MUSTIKA2022}, or \cite{Jo2021}.  Electric vehicles (EVs) are modeled with a vehicle-to-grid function and are therefore modeled similarly to BESS. In addition to the storage constraints, EVs cannot be used when they are not at the charging station. When they return to the charging station, the consumed energy from driving is accounted for, and a minimum SOC is ensured at the next departure time. Finally, the model contains an optimal power flow to represent the grid in which the EC operates. This includes constraints for the power balance and the definition and bounds of the power flow variables. This power flow formulation can be considered standard in this type of literature and can also be found in \cite{WOGRIN2020}. Import and export are modeled as power flows to and from the transmission node. They are bounded by the maximum power of the transformer.  

\subsection{Ex-Post energy allocation}

\begin{figure}
    \centering
    \includegraphics[width=1\linewidth]{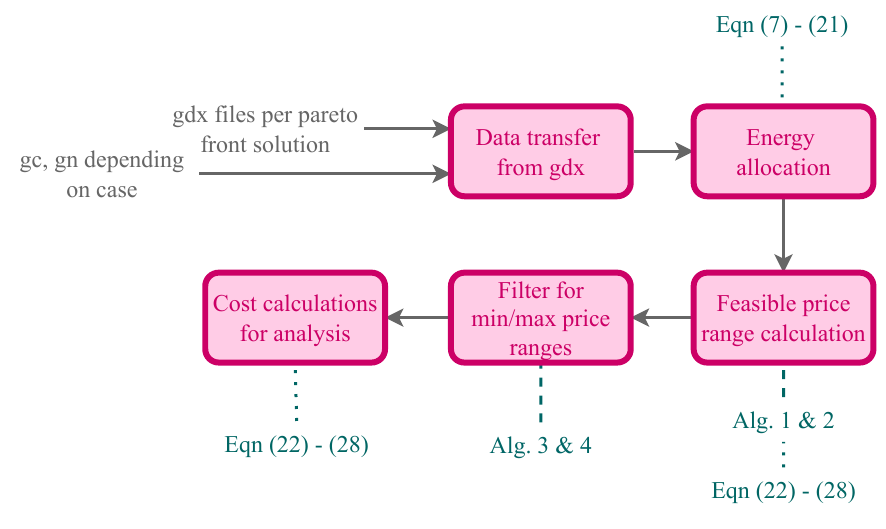}
    \caption{Ex-Post process for analysis of Pareto front impact on the EC and its members.}
    \label{fig:flow}
\end{figure}

The process used to analyze the impact of the Pareto front on the EC and its members is shown in Fig. \ref{fig:flow}. The results from the optimization model come as GAMS-specific gdx-files. Every solution of the Pareto front has its file. Since the model works as a centralized managed EC, the ownership of the individual assets was irrelevant for optimization. It is, however, relevant to the energy allocation and cost calculation. Ownership differs in the cases analyzed in this paper. Therefore, this is an input for the ex-post analysis process. The first step is to transfer the data from the gdx file into a Python commendable pandas data frame. After that, energy allocation between the members is performed using Eqn. (\ref{eqn:NetD})-(\ref{eq:ExeSharen2ec}). 

The energy is allocated according to one of two allocation methods. One is the static method, where a static allocation key is used. The other one is the dynamic method, where the energy is allocated proportionately to each member's share of the overall EC demand. The latter is the most commonly used in Austria and this paper. In the original model, only EVs were member-owned. This paper introduces a new feature that distinguishes between community-owned and member-owned generation and storage in the allocation process. The process, which is done for every timestep h, starts with the calculation of the net demand of each member $D_{h,n}^{net}$ in Eqn. (\ref{eqn:NetD}). In addition to the members' energy demand $D_{h,n}$, consumption $cs_{h,gn}$, and production $p_{h,gn}$ of the generation/flexibility units they own are also accounted for.
When a member has a negative net demand, they act as a producer in this timestep. The total energy supplied by all producing members is summed as $P^{exe,n}$ in Eqn. (\ref{eqn:peev}). This parameter is needed to calculate the EC's net production $P^{c,net}_{h}$ (Eqn. \ref{eqn:Pcnet}), where community-owned consumption and production are considered. The EC may not produce enough energy at a particular timestep to cover the consumption of community-owned storage. A conventional energy supplier then covers this residual grid demand $Rd^c_{h}$ (calculated in Eqn. (\ref{eqn:DcRd})). After this calculation, the net production is set to zero in Eqn. (\ref{eqn:Pcnetzero}) for the timesteps where it is negative. The difference between the total EC demand $D^c_{h}$ (the sum of all members' net demands, Eqn. (\ref{eqn:Dc})) and the EC net production represents the EC excess energy $Exe^c_{h}$, as shown in Eqn. (\ref{eqn:Exed}). The final step before the allocation process involves calculating the actual distributed energy $P^{dis,c}_{h}$, outlined in Eqn. (\ref{eqn:Pdisc}). 
\begin{equation}
\label{eqn:NetD}
    D_{h,n}^{net}=D_{h,n}+\sum_{gn}(cs_{h,gn}-p_{h,gn})~\forall h,n,gn
\end{equation}
\begin{equation}
\label{eqn:peev}
    P_{h}^{exe,n}=\sum_n -D_{h,n}^{net}~\forall h,n, D_{h,n}^{net}\leq 0
\end{equation}
\begin{equation}
\label{eqn:Pcnet}
    P^{c,net}_{h} = \sum_{gc} p_{h,gc} - \sum_{gc}  cs_{h,gc}+P_{h}^{exe,n} ~\forall h,gc
\end{equation}
\begin{align}
%DCBESS
\label{eqn:DcRd}
    Rd^c_{h}=   \begin{cases}
        (-1)\cdot P^{c,net}_{h},                   &for~ P^{c,net}_{h}\leq  0 ~\forall h\\
        0, & \text{else}
                \end{cases}
\end{align}
\begin{equation}
\label{eqn:Pcnetzero}
    P^{c,net}_{h} = 0 ~\forall h,gc, P^{c,net}_{h}\leq  0
\end{equation}
\begin{equation} \label{eqn:Dc}
D^c_{h}=\sum_{n}D_{h,n}^{net}~\forall h,n, D_{h,n}^{net}\geq 0
\end{equation}
\begin{align}
\label{eqn:Exed}
    Exe^c_{h}=   \begin{cases}
        P_{h}^{c,net}- D^c_{h},                   &for~ P_{h}^{c,net}\geq  D^c_{h}, ~\forall h,n\\
        0, & \text{else}
                \end{cases}
\end{align}
\begin{equation}
\label{eqn:Pdisc}
    P^{dis,c}_{h} = P_{h}^{c,net} - Exe^c_{h}~\forall h
\end{equation}
As previously mentioned, the EC energy share of each member is proportional to their share of the EC demand. This is calculated in Eqn. (\ref{eqn:Share,d}), where each member's net demand is divided by the EC demand, yielding a percentage share. The share $Share^{d,\%}_{h,n}$ is then used to determine each member's self-consumption $Sc_{h,n}$, as outlined in Eqn. (\ref{eqn:Sc,d}). Self-consumption refers to the energy each member derives from the EC. At the same time, the residual grid demand $Rd_{h,n}$, as described in Eqn. (\ref{eqn:Rd,d}), represents the energy they still receive from conventional suppliers.
\begin{align} 
\label{eqn:Share,d}
    Share^{d,\%}_{h,n}= \begin{cases}
          \frac{D_{h,n}^{net}}{D^c_{h}}, &for~ D_{h,n}^{net} \geq 0,D^c_{h} \neq 0 ~\forall h,n \\
          0, & \text{else}
    \end{cases}
\end{align}

%\begin{equation} 
%\label{eqn:Sharekw,d}
%Share^{kW}_{h,n}=Shared^{d,\%}_{h,n}\cdot P_{h}^{c,net}%~\forall h,n
%\end{equation}

\begin{equation} 
\label{eqn:Sc,d}
Sc_{h,n}=P^{dis,c}_{h}\cdot Share^{d,\%}_{h,n}~\forall h,n
\end{equation}
\begin{align} 
\label{eqn:Rd,d}
    Rd_{h,n}=\begin{cases}
     D_{h,n}^{net}-Sc_{h,n},& for~D_{h,n}^{net} \geq Sc_{h,n}  ~\forall h,n \\
     0, & \text{else}
    \end{cases}
\end{align}
Suppose there is excess energy in the EC after the energy allocation. In that case, it also needs to be allocated back to a producer so they can sell it to an energy supplier. This is done in a similar manner as the allocation of the energy consumed in the EC. The share of excess energy attributed to each producer corresponds to their share of the total EC energy produced $P_h^c$ (Eqn. (\ref{eq:PC})). The excess energy share of the EC-owned production units $Exe_{h}^{gc}$ is shown in Eqn. (\ref{eq:ExeShareg}) and the one on each member $Exe_{h,n}^{out}$ in Eqn. (\ref{eq:ExeSharen}). The actual energy sold to the EC by each producer $Exe_{h,n}^{in}$ is calculated in Eqn. (\ref{eq:ExeSharen2ec}) by subtracting the excess energy sold to their supplier from their original surplus energy presented by their net demand. 
\begin{equation}
    \label{eq:PC}
    P_{h}^{c} = \sum_{gc} p_{gc} + P_{h}^{exe,n}
\end{equation}

\begin{align} 
    \label{eq:ExeShareg}
    Exe_{h}^{gc}= \begin{cases}
    Exe_{h}^{c} \cdot \frac{\sum_{gc} p_{gc}}{P_{h}^{c}} &for~P_{h}^{c} \neq 0,~\forall h,gc \\
     0, & \text{else}
    \end{cases}
\end{align}
\begin{align} 
    \label{eq:ExeSharen}
    Exe_{h,n}^{out}= \begin{cases}
    Exe_{h}^{c} \cdot \frac{(-1)\cdot D_{h,n}^{net}}{P_{h}^{c}} &for~P_{h}^{c} \neq 0,D_{h,n}^{net} \leq 0~\forall h,n \\
     0, & \text{else}
    \end{cases}
\end{align}
\begin{equation}
\label{eq:ExeSharen2ec}
    Exe_{h,n}^{in}=(-1)\cdot D_{h,n}^{net}-Exe_{h,n}^{out}~\forall h,n,~D_{h,n}^{net} \leq 0
\end{equation}

\subsection{Ex-Post cost calculation}

Finally, the energy is allocated, the results are used to determine the cost of each member and the EC. Here, $Cost_n^{in}$ is the expenditure from energy in the EC and $Cost_n^{out}$ for energy from a supplier. The members pay the overhead costs for metering for the number of metering points they possess. For the energy bought from the EC, they pay a price set by the EC as a collective and grid tariffs (see Eqn. (\ref{eqn:Costinn})). If they sell energy to the EC, they will get their excess reimbursed with $C^{EC,g}$. Usually, the consumers pay less than the producers earn. The difference gives the EC funds to pay for the remaining overhead costs. The cost of energy paid to the supplier is determined in Eqn. (\ref{eqn:Costoutn}) includes the consumption-independent component of the grid tariffs, the cost from the residual grid demand, including the corresponding grid tariffs, and the profit from selling the excess energy. 

\begin{align}
\label{eqn:Costinn}
    Cost^{in}_{n}=C^{ovh,met}(1+card(gn)) +\sum_h \Big( \\ \nonumber
    Sc_{h,n}\cdot (C^{EC,n} + C^{gt,in})-Exe_{h,n}^{in}\cdot C^{EC,g}\Big)~,\forall h,n
\end{align}
\begin{align}
\label{eqn:Costoutn}
    Cost^{out}_{n}= C^{gt,flat}+
    \sum_h \Big( Rd_{h,n}\cdot ( C^{import}_h+ C^{gt,out}) \\ \nonumber
    - Exe_{h,n}^{out}\cdot C^{export} \Big)~,\forall h,n
\end{align}
The cost of the community as an entity is very similar to that of the members and can be differentiated between costs inside $Cost^{c,in}$ and outside the EC $Cost^{c,out}$. As the EC-owned assets work for the community as a whole, they are not affected by the EC energy price. However, for energy consumed by the community-owned BESS grid tariffs must be paid. Again, there is a difference between grid tariffs for energy from the community and the supplier. If there is a surplus of energy from community-owned BESS or PV, it is sold to the supplier. The third type of cost that occurs for the EC is the overhead cost presented in Eqn. (\ref{eqn:costcov}). As already described earlier, they are set up as fixed accounting costs, variable overhead costs depending on the amounts of allocated energy in the EC, and costs per metering point. The just elaborated costs of the EC need to be offset by the difference between producer and consumer energy prices. This balance is shown in Eqn. (\ref{eqn:costcbalance}). If the EC decides to subsidize its energy-poor members, they do not pay the EC energy price and no metering point overhead cost. In this scenario, the balance must be attained only with the cost contribution from the other members (Eqn.(\ref{eqn:costcbalancepoor})). 
\begin{equation}
\label{eqn:costcin}
    Cost^{c,in}= \sum_{h} \Big( \sum_{gc} cs_{h,gc} -Rd_{h}^{c} \Big)\cdot (C^{gt,in})~,\forall h
\end{equation}

\begin{align}
\label{eqn:costcout}
    Cost^{c,out}=  C^{gt,flat}+ \\ \nonumber \sum_h \Big( Rd_{h}^{c}\cdot ( C^{import}_h+ C^{gt,out}) 
    - Exe_{h}^{gc}\cdot C^{export} \Big)~,\forall h
\end{align}

\begin{align}
\label{eqn:costcov}
    Cost^{c,ovh}=  C^{ovh,acc}+  C^{ovh,var}\cdot\sum_{h,c}Sc_{h,n}\\ \nonumber
    + C^{ovh,met}\cdot card(n) \cdot card(g)~,\forall h
\end{align}
\begin{align}
\label{eqn:costcbalance}
    Cost^{c,net}=  Cost^{c,in}+  Cost^{c,out} + Cost^{c,ovh}\\ \nonumber
    -\sum_{h,n}\Big(Sc_{h,n}\cdot C^{EC,n}- Exe_{h,n}^{in}\cdot C^{EC,g}\Big) \\ \nonumber
    - C^{ovh,met}\cdot card(n) \cdot card(g)~,\forall h,n
\end{align}
\begin{align}
\label{eqn:costcbalancepoor}
    Cost^{c,net,ep}=  Cost^{c,net}~,\forall h,n \notin poverty
\end{align}
\subsection{Ex-post feasible EC energy price range analysis}

\begin{figure}
    \centering
    \includegraphics[width=0.6\linewidth]{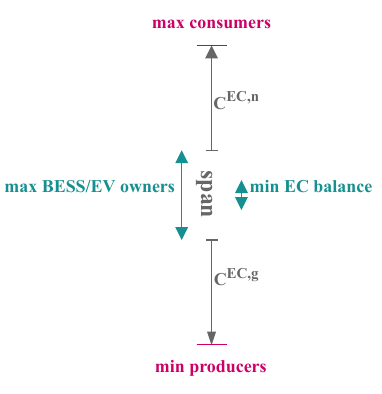}
    \caption{EC energy price range boundaries.}
    \label{fig:price_range}
\end{figure}

Price setting is not trivial in an EC. Other than the fact that this process should be done as a collective, everyone has their interests. In finding a feasible range for $C^{EC,n}$ and $C^{EC,g}$, this paper defines prices as feasible if the EC cost balance is positive and every member has a better cost outcome than not being in an EC. How this affects the setting of the feasible price range is illustrated in Fig. \ref{fig:price_range}. As shown in the figure, the maximum $C^{EC,n}$ is determined by the consumers, and the minimum $C^{EC,g}$ by the producers. As mentioned, the price difference needs to be at least enough for the EC cost balance to be positive but not too far apart for the BESS and EV owners to be worse off than not being in an EC. In order to identify the feasible price ranges, Algorithms 1 and 2 were created. With them, the six highest and six lowest feasible ranges were located. Algorithms 3 and 4 were used to find the ranges with maximum and minimum price spans of the ones found in the first two algorithms. This results in one upper and one lower price range for each Pareto front solution. The pseudo code for all algorithms can be found in the Appendix. 
 
\section{Case Studies}
\label{sec:CS}

\begin{figure}
    \centering
    \includegraphics[width=1\linewidth]{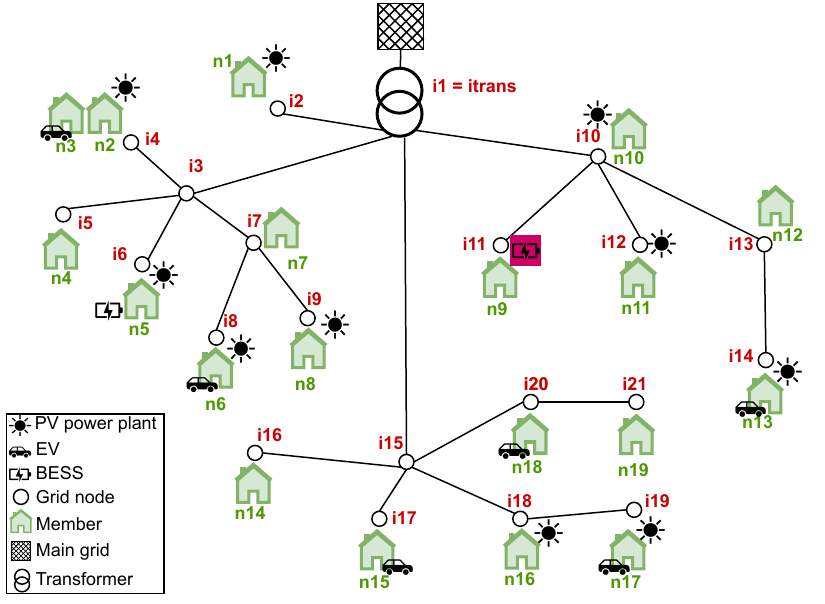}
    \caption{Case study EC including grid topology, PV, EV and BESS allocation do members.}
    \label{fig:grid}
\end{figure}

This section presents the stylized EC design chosen to implement the methodology above. The main characteristics are described here, while the exact input parameters can be found in the online appendix \cite{github}. All analyzed scenarios are based on a stylized local REC illustrated in Fig. \ref{fig:grid}. It has 19 household members with six EVs, ten PV units, and two BESSs. The low-voltage grid model that the EC is set up in is a 21-bus system loosely based on a partial distribution grid in western Austria. Bus 1 is the transmission node where energy is imported and exported to the main grid. In reality, this would be the low-voltage side of the transformer. Smart meter data from the Open Power System Data project \cite{OPSD2020} were implemented to model the households. Because of the limited number of demand profiles in the source, they were scaled by a factor to create 19 unique demand patterns. The EC has four 10.22~kWp and six 5.11~kWp PV units. In order to ensure a closer correlation between demand and production profiles, capacity factors from south Germany, where the smart meter data is from, were used. The PV production data was obtained from  Renewables Ninja \cite{Pfenninger_2016}. The BESS units in the model have a capacity of 11~kWh with a charging efficiency of 96\%. The EVs in the EC were modeled after the four most popular EVs in Austria, according to Statistics Austria \cite{EVstat}. The member's departure and arrival times differ per household and lie between 6:00–9:00 and 15:00–19:00. The energy the EVs consume is calculated with travel distance data from Probst et al. \cite{Probst2011}. The grid tariffs are from one of Austria's larger DSOs' price sheets. Import and Export prices were derived from an Austrian supplier \cite{awttar}. Dynamic rates were chosen to provide price incentives. For importing energy, this is hourly, and for exporting monthly, because at this point, no supplier offers hourly rates for PV energy. Lastly, the overhead costs are provided by a small Austrian EC service provider.

Using this proposed EC constellation, three case studies were analyzed. The first case serves as reference where all PVs, EVs, and BESS are owned by members. In a second case study, the value of a community owned BESS is explored. For this instance, the BESS at node $i11$ (highlighted in Fig. \ref{fig:grid}) is community-owned instead of member $n9$. The last case study subsidizes two energy-poor members ($n7$ and $n14$) by exempting them from the EC energy price and the metering point overhead costs. Those members only pay the reduced grid tariff for the energy they consume from the EC.  

\section{Results}
\label{sec:Results}

\subsection{Impact on EC operation}

\begin{figure}
    \centering
    \includegraphics[width=0.7\linewidth]{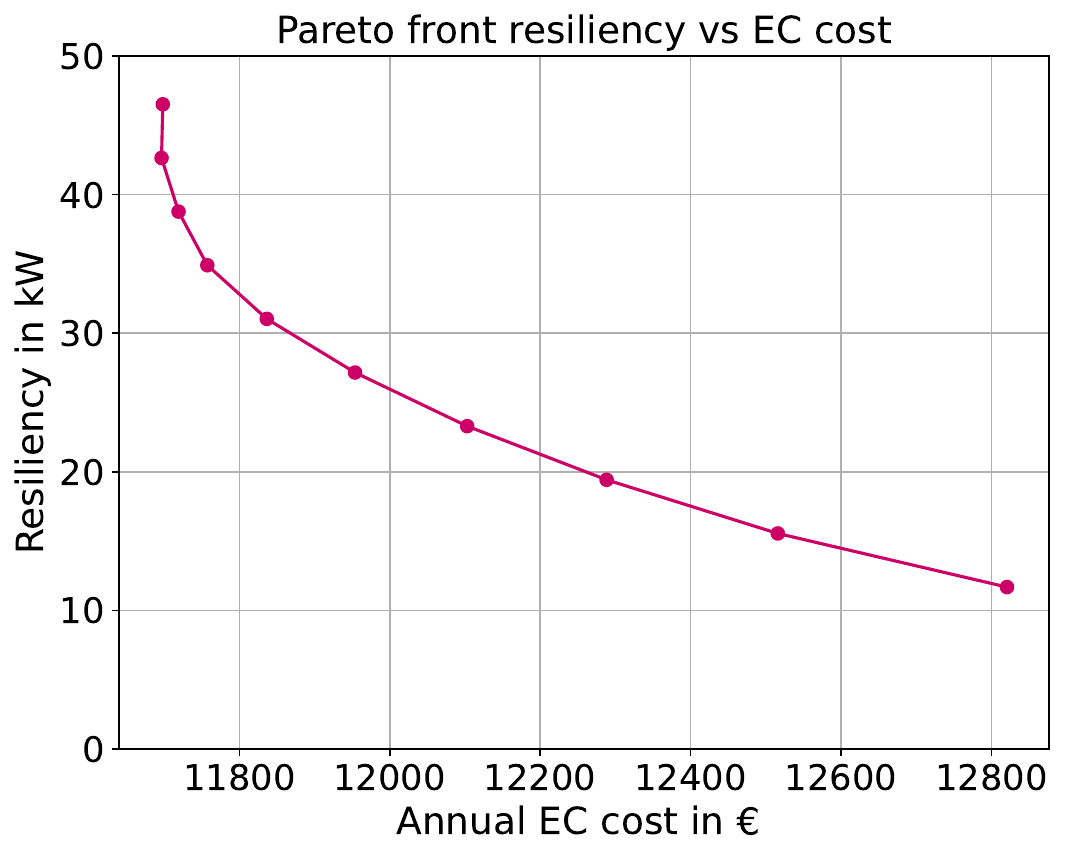}
    \caption{Pareto front of EC cost and resilience.}
    \label{fig:pareto}
\end{figure}

In Fig. \ref{fig:pareto}, the results of the multi-objective optimization are illustrated by the Pareto front between the annual EC cost and resilience. The EC doesn't have enough flexibility to be completely independent of the main grid. 
However, maintaining such independence can be quite costly for the EC. While a smaller impact on the higher grid levels is favorable for the DSO, there is no benefit from EC's (theoretical) complete autarky. The primary goal is to avoid creating significant power peaks in either direction, as these must be addressed at the higher grid levels. The Pareto front shows that the additional financial effort increases with every solution closer to the minimum power peak. This suggests that small financial incentives would make it worthwhile for the EC to adjust its operations accordingly. 

\begin{figure}
    \centering
    \includegraphics[width=\linewidth]{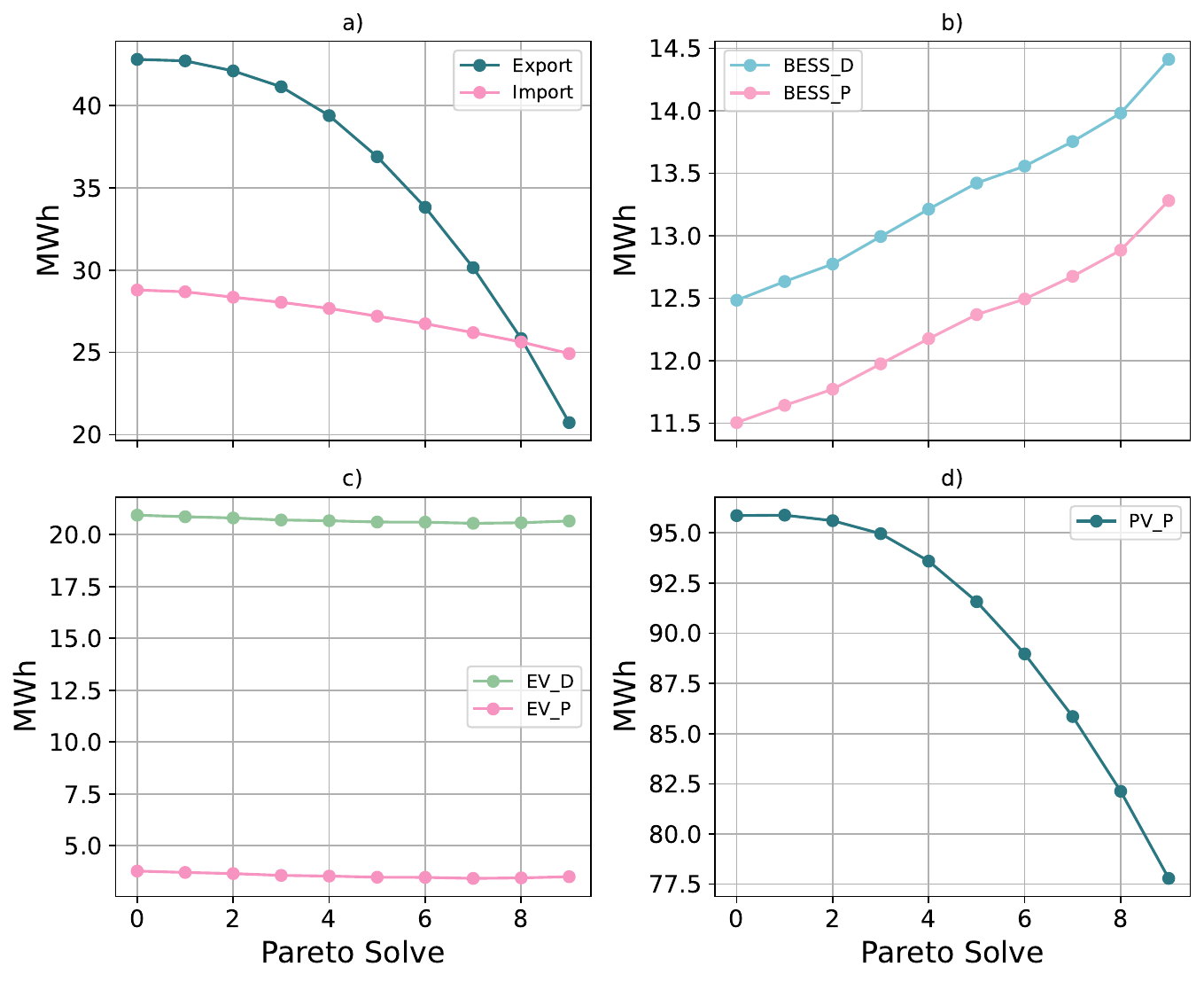}
    \caption{Operational Pareto front results: a) imports and exports, b) BESS production and demand, c) EV production and demand, d).}
    \label{fig:sumpareto}
\end{figure}

Fig. \ref{fig:sumpareto} shows the impacts on production and consumption units. As expected, as resilience improves, annual imported and exported energy decreases. The sharper decrease of the export can be explained with the curtailment option present in the model. While demand needs to be met (and, if necessary, imported), surplus power can be curtailed "only," resulting in a loss of profit. This is also represented by the PV production curve, which is almost identical in shape to the export curve. Another expected outcome is an increase in BESS utilization—both demand and production of the BESSs in the EC increase with improved resilience. However, EV production and demand remain essentially unchanged. The primary cause for this is their unavailability during peak PV production times.
\begin{figure}
    \centering
    \includegraphics[width=1\linewidth]{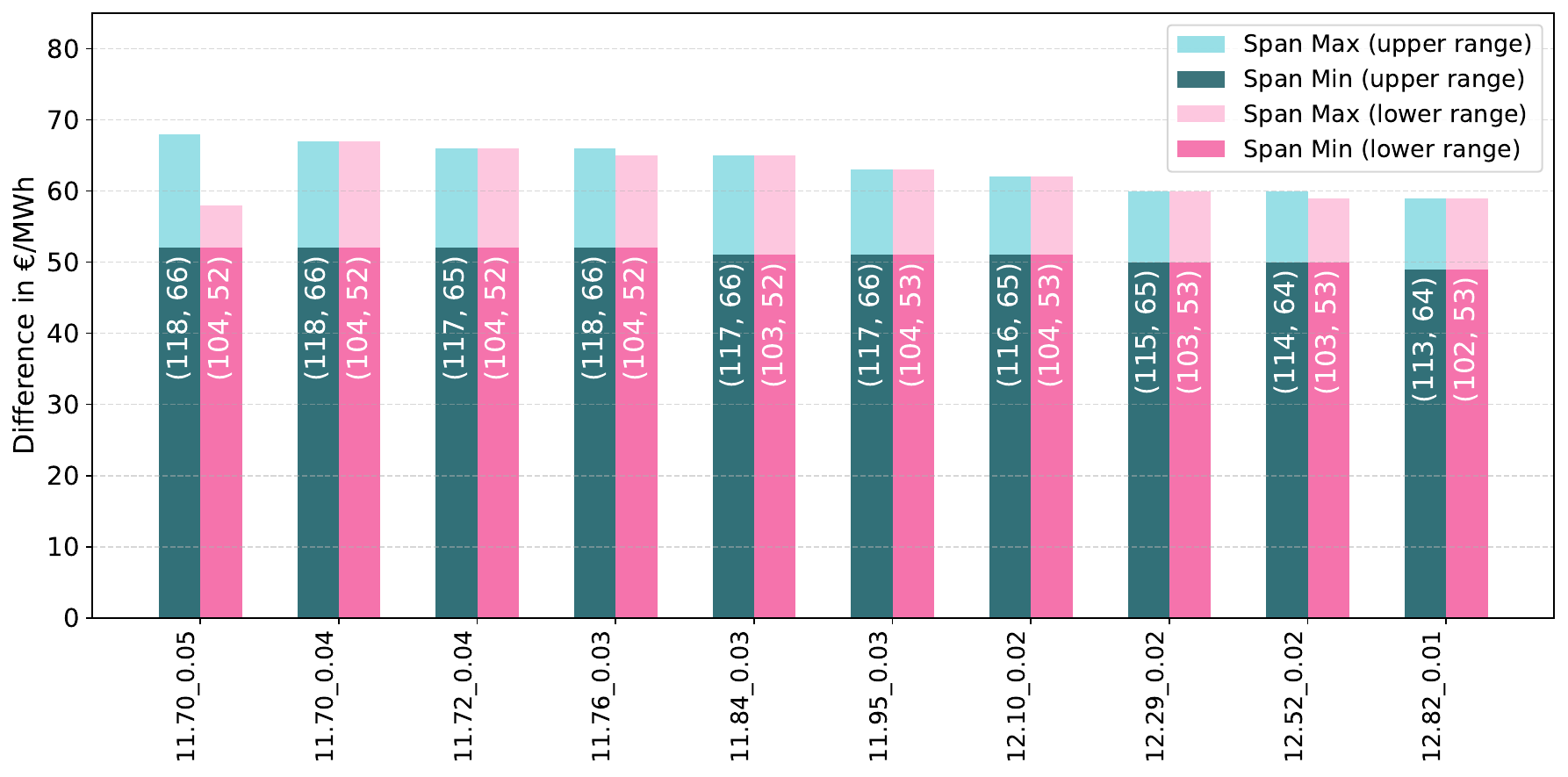}
    \caption{EC energy price spans between consumer and producer price for the reference case.}
    \label{fig:CS4_range}
\end{figure}
\subsection{Impact on feasible price ranges}
\begin{figure}
    \centering
    \includegraphics[width=1\linewidth]{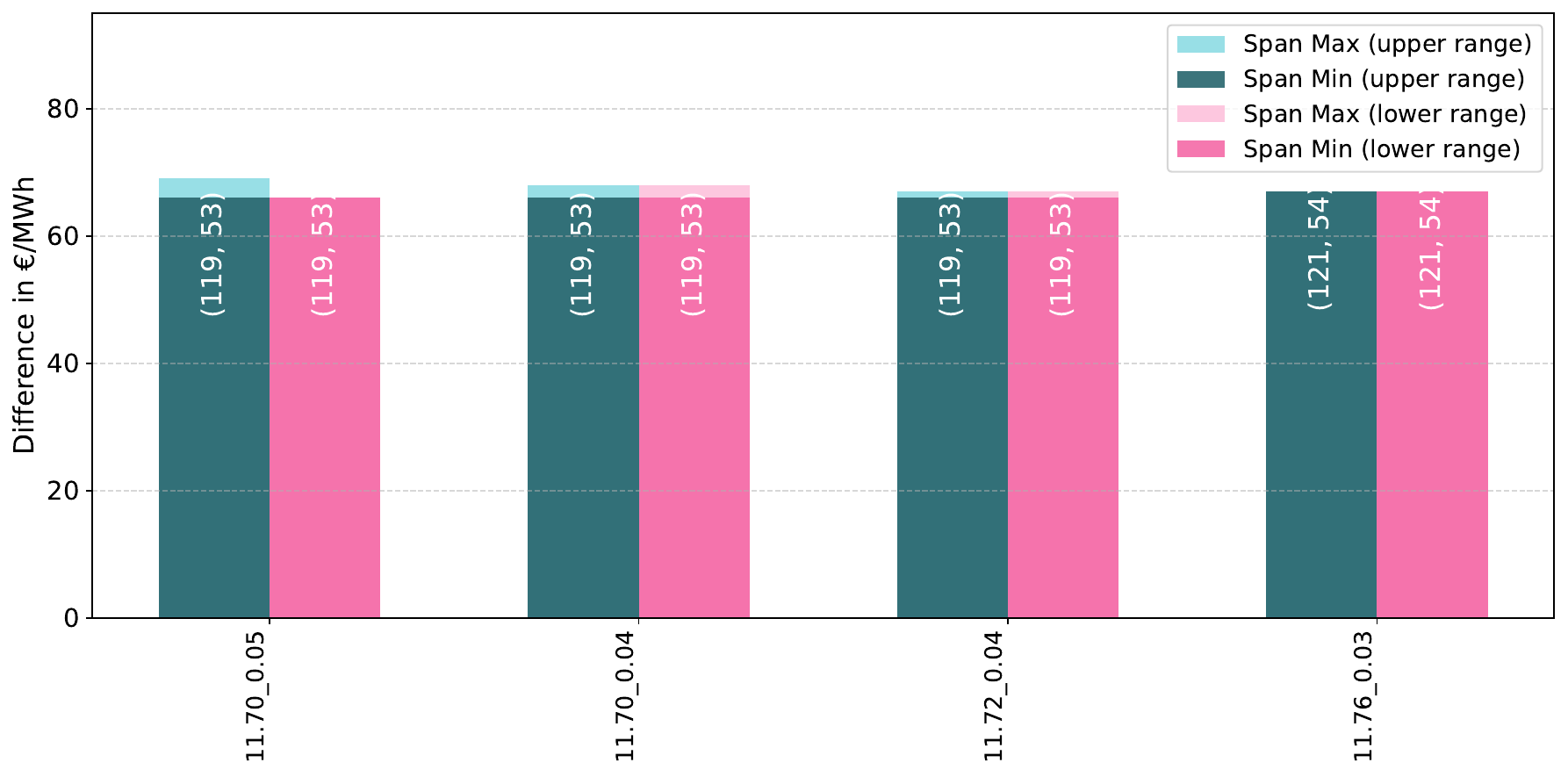}
    \caption{EC energy price spans between consumer and producer price for the communal BESS case.}
    \label{fig:CS3_range}
\end{figure}

\begin{figure}
    \centering
    \includegraphics[width=1\linewidth]{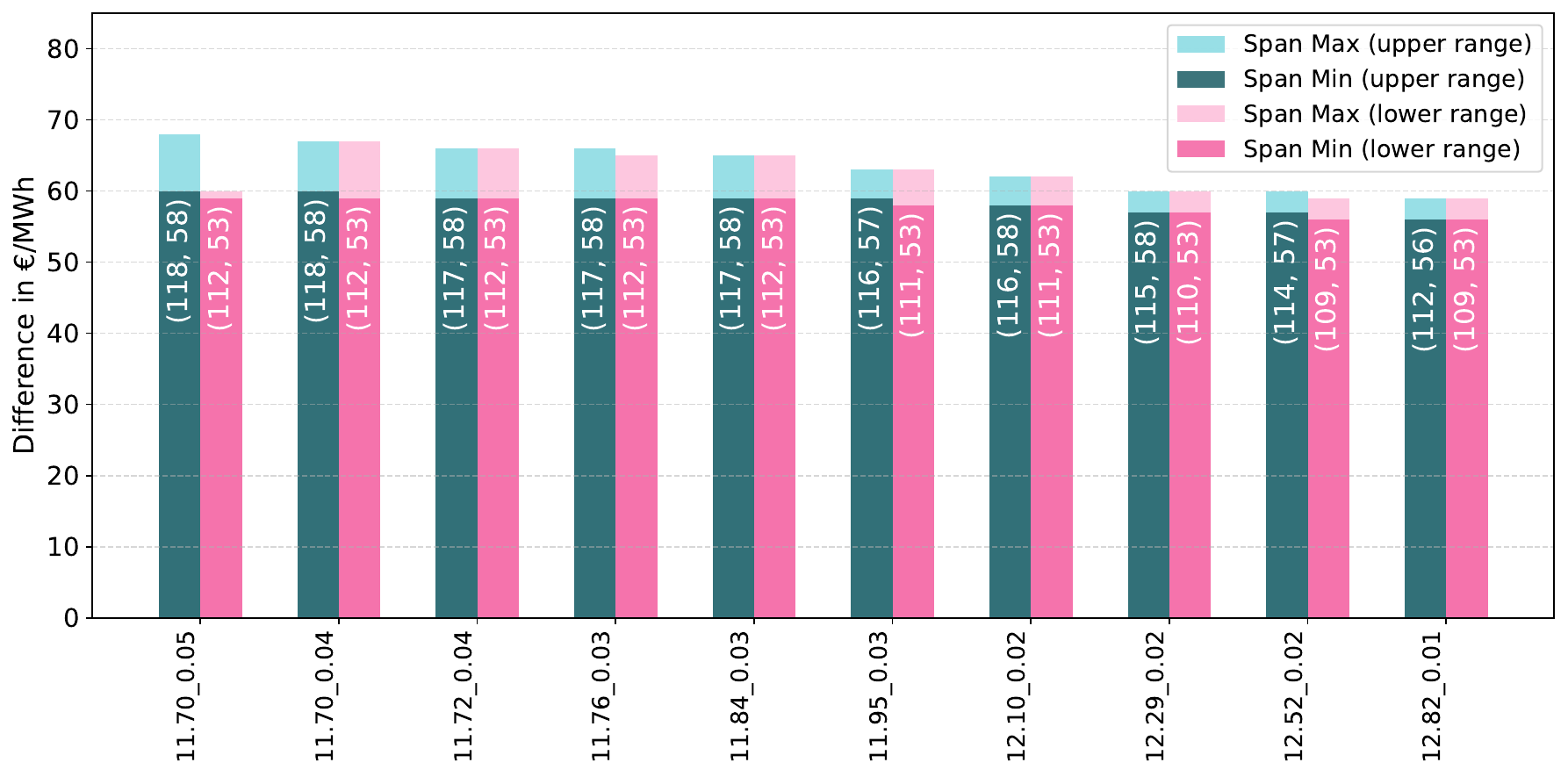}
    \caption{EC energy price spans between consumer and producer price for the energy poverty subsidy case.}
    \label{fig:CS5_range}
\end{figure}

The impact of the trade-off between resilience and cost can be observed in the feasible ranges for energy prices related to buying and selling in the EC. The results indicate that when moving along the Pareto front, the degrees of freedom in setting the EC energy prices sink. The maximum and minimum span between consumption price $C^{EC,n}$ and production price $C^{EC,g}$ both decrease with more resilient solutions. This effect is more pronounced in the maximum span than in the minimum span. These spans, including the corresponding minimal $C^{EC,n}$ and $C^{EC,g}$, are shown in Fig. \ref{fig:CS4_range}-\ref{fig:CS5_range} for the three cases. A decline of $C^{EC,nmax}$ can be observed in the upper price range, while a slight increase of $C^{EC,gmin}$ can be observed in the lower price range. Among the three cases, the impact of the Pareto front was the biggest for the communal BESS. Balancing the price into a feasible region was only possible for the four more cost-efficient solutions. As expected the energy poverty subsidy case was more challenging to balance than the reference case. Consequently, this resulted in broader spans with lower $C^{EC,n}$ in the upper price range and higher $C^{EC,g}$ in the lower price range.

\subsection{Impact on member groups}
The maximum price span is helpful if the EC plans to invest in communal production or flexibility assets. Otherwise, the EC does not need to make a lot of profit. Since potential investments are not in the scope of this paper and the minimum span results in a smaller financial burden on the members with EVs and BESS, the analysis focuses solely on the minimum price span. The results are presented both for the lower and upper price ranges. The overall trends of how the Pareto front impacts the cost of the member groups are consistent across the case studies, although the magnitude of these effects varies. The financial benefit ranking of the three cases is also consistent among member groups (except for the energy-poor members) but again varies in magnitude. Looking at this ranking, an interesting finding from this study is that the communal BESS did not provide financial benefits to any member group. One of the reasons for this is the larger price span required for EC energy prices. While the theoretical advantage of a communal BESS is that it allows for more community-shared energy—since member-owned systems only share net production (BESS production minus owner demand)—this benefit was not sufficient to offset the necessary pricing for a feasible price range, leading to an overall negative impact. The third case of subsidizing energy poor members does impact all other member groups, compared to the reference case. 

\begin{figure}
    \centering    \includegraphics[width=1\linewidth]{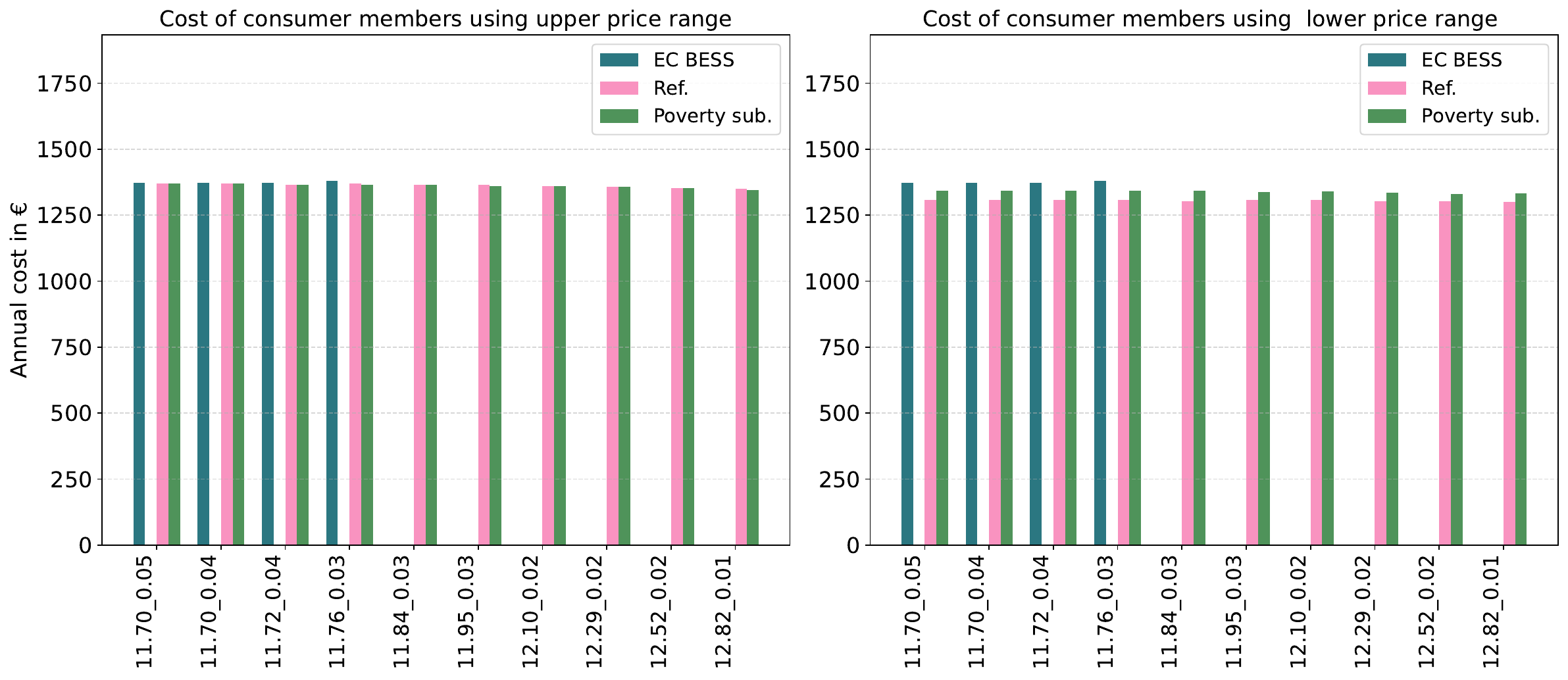}
    \caption{Summarized cost for consumer members.}
    \label{fig:CS4_consumers}
\end{figure}

\begin{figure}
    \centering    \includegraphics[width=1\linewidth]{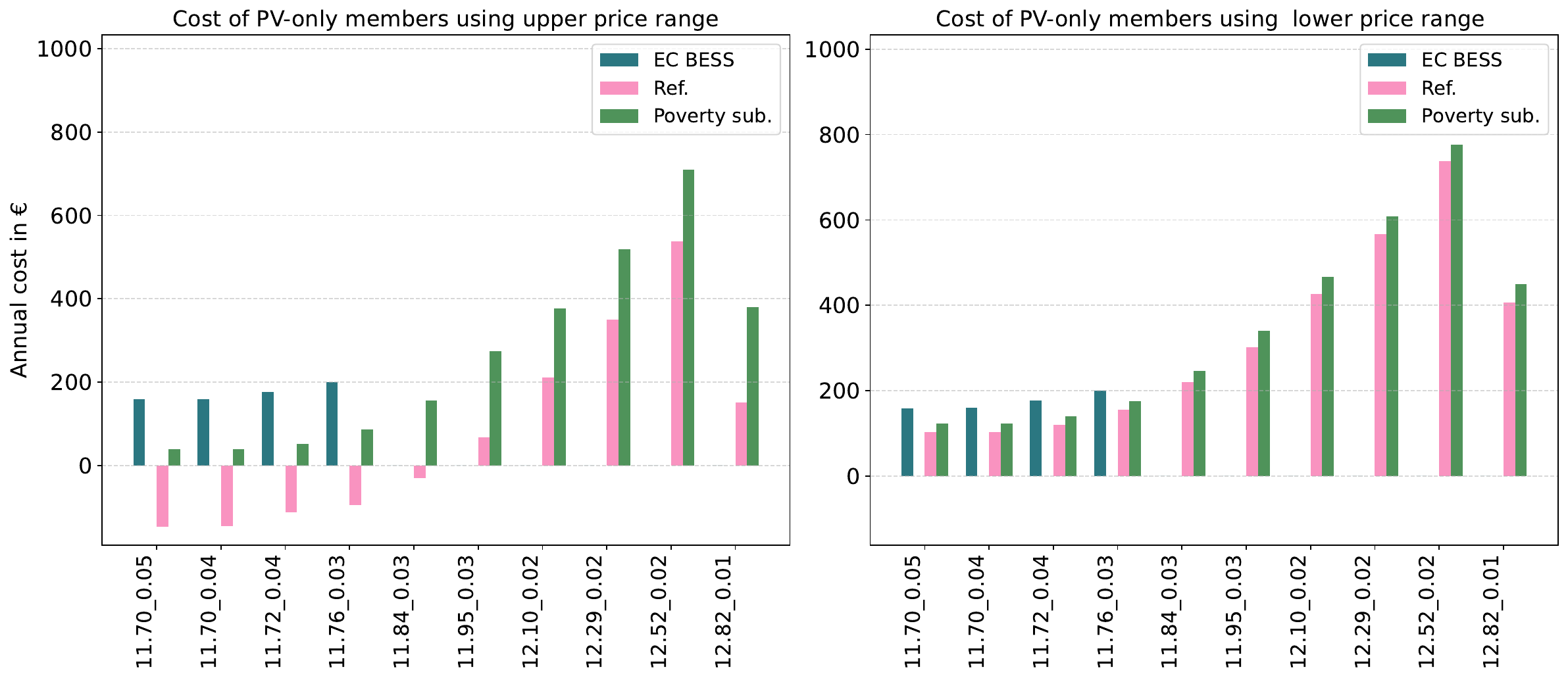}
    \caption{Summarized cost for PV-only members.}
    \label{fig:CS_PVonly}
\end{figure}

To ensure consistency when analyzing the member groups, member $n9$ is excluded from the plots because they are a consumer in two cases and a BESS owner in the other. Also, $n7$ and $n14$ were excluded from the consumer group and analyzed as a separate group of energy-poor members. Fig. \ref{fig:CS4_consumers} shows consumers are largely unaffected by moving along the Pareto front. Additionally, there is minimal difference among the case studies in the upper price range. In the lower price range, however, the negative impact of the communal BESS increases slightly. As expected, the energy poverty case is more expensive than the reference case, but not by a large margin. Fig. \ref{fig:CS_PVonly} illustrates the cost of the members that only own a PV system. Members with a combination of EVs or BESSs will be analyzed in their own category. For PV owners, the trade-off for increased resilience has a negative impact on their annual costs. Due to curtailed energy, they save less money by selling excess energy to the supplier.

\begin{figure}
    \centering    \includegraphics[width=1\linewidth]{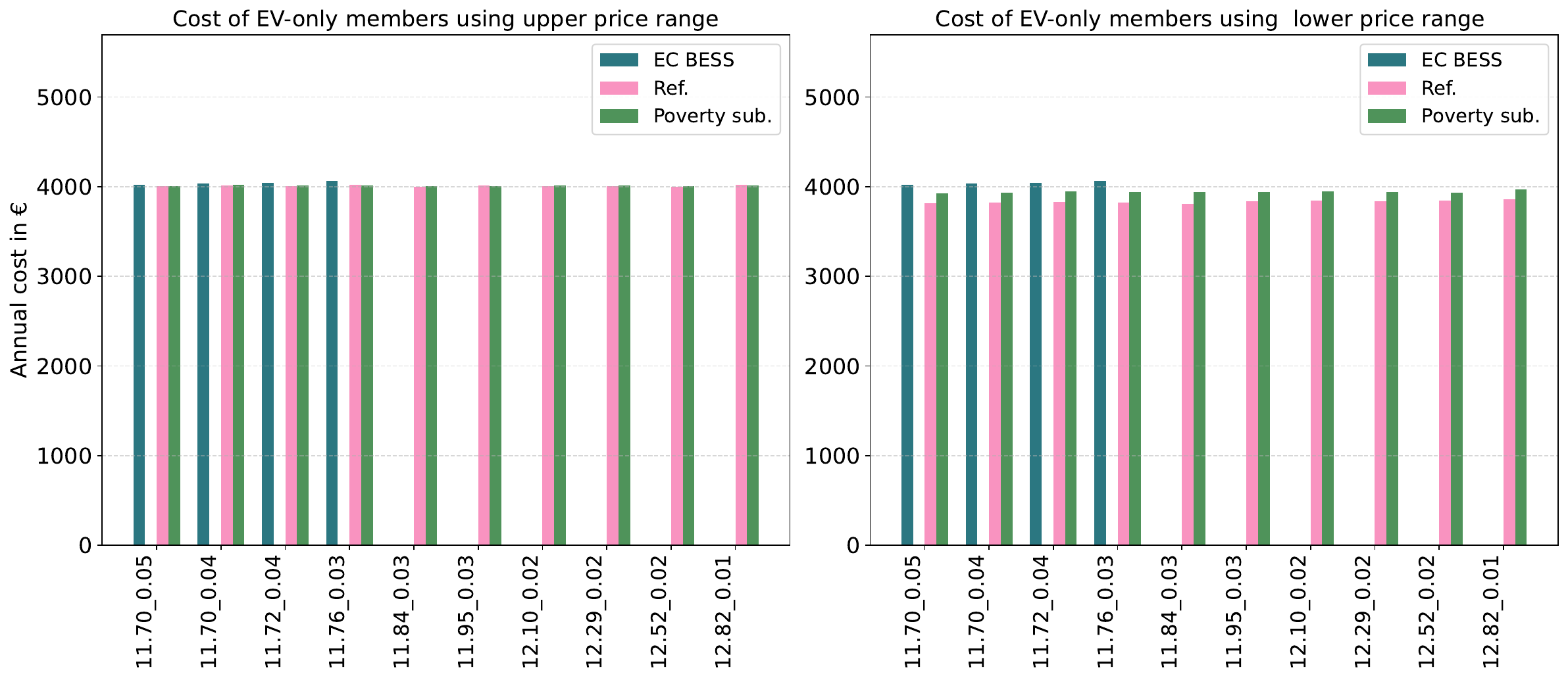}
    \caption{Summarized cost for EV-only members.}
    \label{fig:EVonly}
\end{figure}

\begin{figure}
    \centering    \includegraphics[width=1\linewidth]{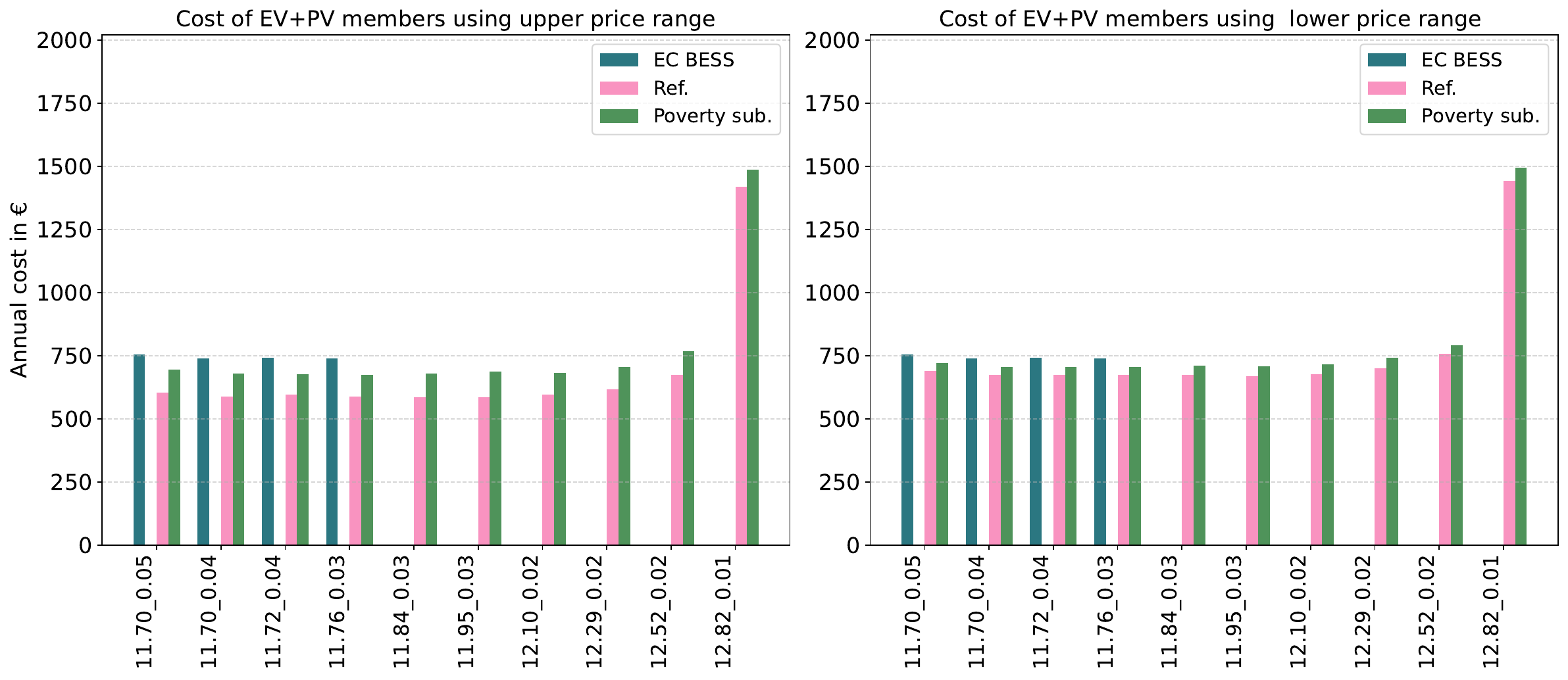}
    \caption{Summarized cost for PV and EV owning members.}
    \label{fig:CS_EVPV}
\end{figure}

Nevertheless, in the reference case's upper price range, they still profit from the more cost-minimized half of the Pareto front. It is also worth noting that the difference between the energy poverty subsidy case and the reference is not insignificant in the upper price range. The previous results indicated that EVs were not heavily used for flexibility. Therefore, the EV-owner member group results are similar to mere consumer members, as illustrated in Fig. \ref{fig:EVonly}. Members that own a PV unit in addition to their EVs have a relatively stable result along the Pareto front, except at the maximum resilience solution, as depicted in Fig. \ref{fig:CS_EVPV}. For this member group, the difference between cases is more pronounced in the upper price range in comparison to the EV owners without a PV unit.


\begin{figure}
    \centering    \includegraphics[width=1\linewidth]{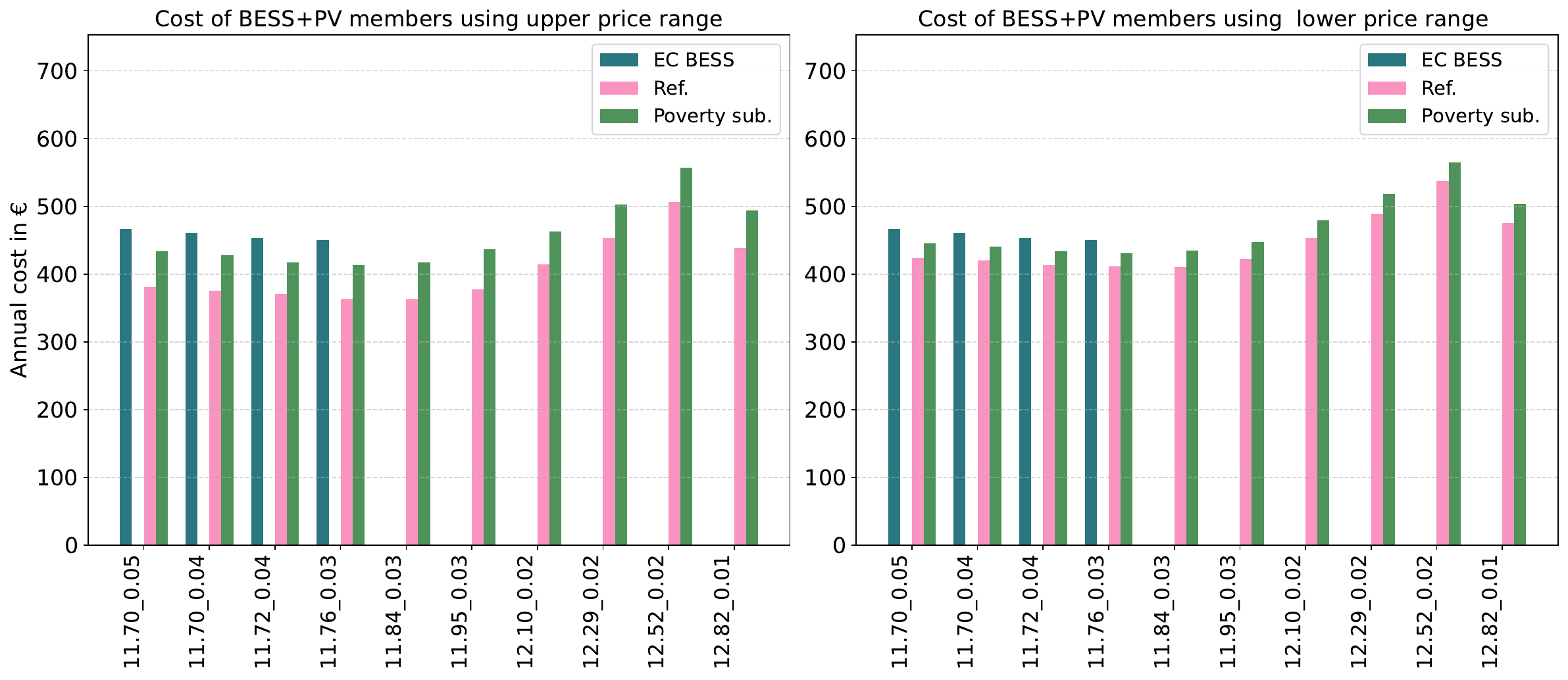}
    \caption{Summarized cost for PV and BESS owning members.}
    \label{fig:CS_BESSPV}
\end{figure}

Due to the more flexible use of BESS in comparison to EVs, the results of members with a combination of PV and BESS follow a more similar trend to the PV owners, which is shown in Fig. \ref{fig:CS_BESSPV}. In this case, costs increase with greater resilience, like for members with only PV, but less severe. Another interesting result is that the curtailment on the resilience extreme of the Pareto front is not distributed homogeneously between members. While overall PV production decreases, generation actually increases for most PV-only members. Conversely, members with both PV and EVs experience the highest levels of curtailment. This observation explains the cost trend reversal on the last Pareto front solution in Fig. \ref{fig:CS_PVonly}, \ref{fig:CS_EVPV}, and \ref{fig:CS_BESSPV}. 

\begin{figure}
    \centering    \includegraphics[width=1\linewidth]{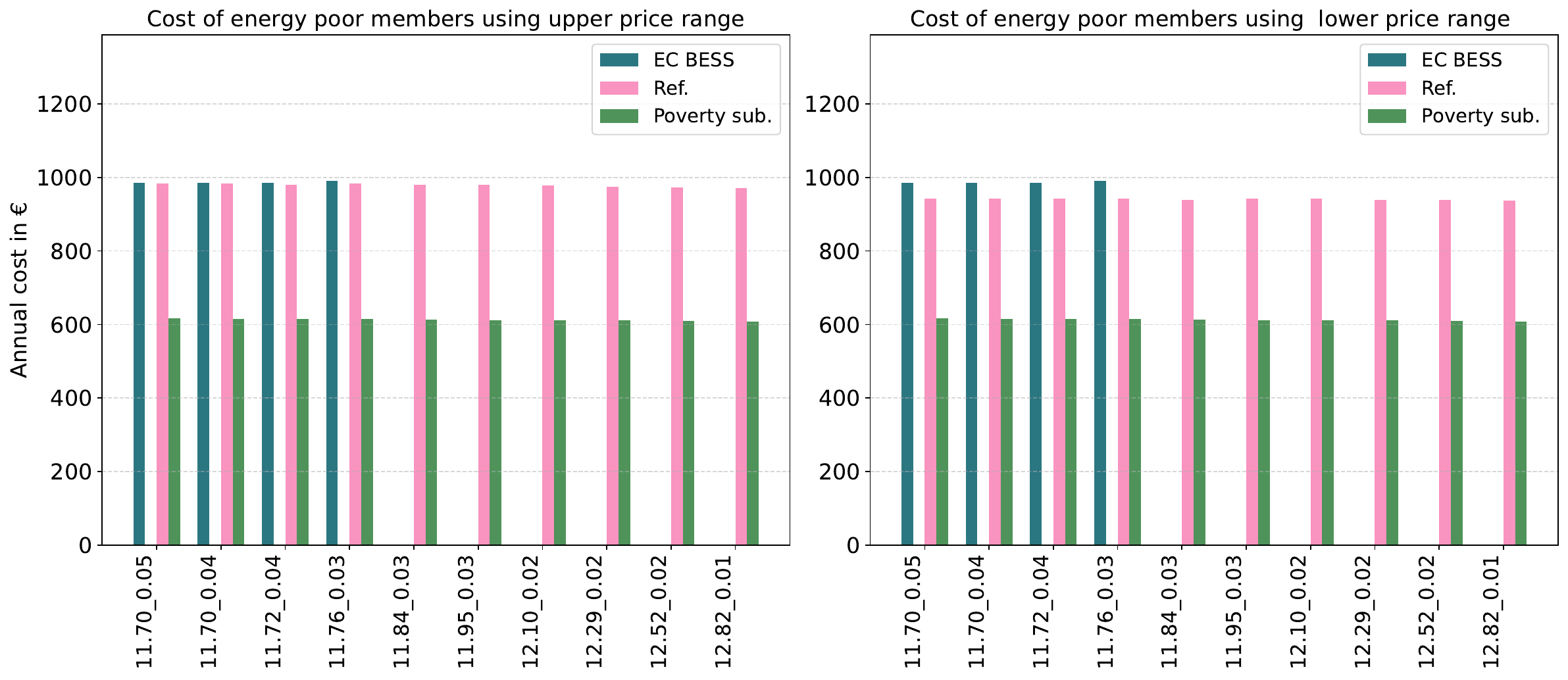}
    \caption{Summarized cost for energy poor members.}
    \label{fig:CS_Poor}
\end{figure}

Lastly, the total costs for this group are illustrated in Fig. \ref{fig:CS_Poor} to analyze the effectiveness of the subsidies provided to energy-poor members. For $n7$, a cost decrease of approximately 39\% in the upper price range and 36\% in the lower price range compared to the reference case was achieved. The results for $n14$ were similar, with 35\% and 33\%, respectively.

\section{Conclusion}
\label{sec:conclusion}

This paper developed a multi-objective optimization framework to explore the trade-offs between economic and technical objectives in ECs. Using the augmented $\epsilon$-constraint method, the annual total cost of the EC and its grid resilience, represented by the reduction of power peaks at the transformer, were optimized. While more increased resilience could be reached, it comes at a higher operational cost for the EC. The Pareto frontier analysis suggested that initial improvements in resilience can be achieved with relatively modest cost increases from the EC. However, as resilience continues to increase, the additional costs needed for incremental gains in resilience rise significantly. A deeper analysis of the results showed that the financial repercussions of resilience enhancement are unevenly distributed among member groups. Members who own PV panels are financially impacted by increased resilience. PV production is affected through curtailment, leading to a profit loss. Mere consumers, however, are almost unaffected. The trends of how the trade-off affected the establishment of feasible EC energy prices stayed the same for all analyzed EC operation strategies. The intensity, however, differed, making implementing a communal BESS at high levels of resilience not financially feasible in the analyzed EC. The communal BESS also did not result in direct financial benefits for individual members under the energy allocation and cost mechanism used in this paper. Also, subsidizing energy-poor members made finding feasible price ranges more challenging than in the reference case. 

Nevertheless, the subsidy provided the energy-poor members with a decrease in overall energy cost of more than 30\%. From a policy perspective, the findings suggest that encouraging grid-friendly EC operations could improve resilience at the distribution level. The results also underscore the importance of analyzing the more granular effects inside the EC instead of only considering the EC as one unit. This paper also concludes that EC policies for utilizing their technological potential should not come at the expense of vulnerable households. Future research will attempt to include energy poverty abatement as an objective function. Implementing this into the multi-objective optimization model brings the research from only analyzing the impact on energy-poor members to including it in EC operational decision-making.  
%\section*{Acknowledgments}

\bibliographystyle{IEEEtran}
 \bibliography{MPiEC_IEEE_Gruber}
%\section*{Biography Section}
%\bf{If you include a photo:}\vspace{-33pt}
%\begin{IEEEbiography}[{\includegraphics[width=1in,height=1.25in,clip,keepaspectratio]{fig1}}]{Michael Shell}
%Use $\backslash${\tt{begin\{IEEEbiography\}}} and then for the 1st argument use $\backslash${\tt{includegraphics}} to declare and link the author photo.
%Use the author name as the 3rd argument followed by the biography text.
%\end{IEEEbiography}

\section*{Data availability}
\noindent As stated in the paper, the data used in the case study can be found in an online appendix:\\
\textcolor{blue}{https://github.com/liaGTUG/Cost-versus-Resilience-in-Energy}\\
\textcolor{blue}{-Communities-A-Multi-Objective-Member-Focused-Analysis}
\section*{Appendix}
\begin{comment}
\begin{algorithm}[H]
\caption{Feasible upper price range}
\State $mode \gets on$
\State $counter \gets off$
\State $cn \gets c_{n,input}$
\State $cg \gets c_{g,input}$
\While{$cg > c_{g,min}$}
    $cg = cn$ 
    \While{$cg > c_{g,min}$}
        compute cost
        \If{$[cost_{n,with EC} < cost_{n,without EC}~\forall n]~\&~[cost_c < 0]$} 
            \If{$mode = = off$}
                safe $(cg,cn)$
                $mode = on$
                \If{$counter = = 0$}
                    \State $counter \gets 1$
                \EndIf
            \EndIf
            \State $cg \gets cg-1$   \\
            $i \gets i-1$
        \Else
            \If{$mode == off$}
                \State $i \gets i+2$
            \EndIf
        \EndIf
    \EndWhile
\EndWhile    
\end{algorithm}
\end{comment}
\begin{algorithm}[H]
\caption{Feasible upper range}
\begin{algorithmic}[1]
    \For{solv \textbf{in} Pareto}
        \State mode $\gets$ on
        \State counter $\gets$ 0
        \State $cn \gets cn^{input}$
        \State $cg^{min} \gets cg^{input}$
        
        \While{$cn > cg^{min}$}
            \State $cg \gets cn$ 
            \While{$cg > cg^{min}$}
                \State \textbf{compute} $COST$
                \If{\( [cost_{n}^{wEC} < cost_{n}^{woEC}~\forall n]~\land~[cost_c < 0] \)}
                    \If{mode = off}
                        \State mode $\gets$ on
                        \State \textbf{save} $f_{solv}^{up}(cn_{counter}\gets cn)$
                        \State \textbf{save} $f_{solv}^{up}(cg_{counter,mode} \gets cg)$
                        \If{counter = 0}
                            \State counter $\gets$ 1
                        \EndIf
                    \EndIf
                    \State $cg \gets cg - 1$
                \Else
                    \If{mode = on}
                        \State mode $\gets$ off
                        \State \textbf{save} $f_{solv}^{up}(cg_{counter,mode} \gets cg+1)$
                        \State \textbf{break}
                    \EndIf
                    \State $cg \gets cg - 1$
                \EndIf
            \EndWhile
            \State mode $\gets$ off
            \State $cn \gets c_{n}-1$
            \If{counter $> 0$}
                \State counter $\gets$ counter $+1$
            \EndIf
            \If{counter = 6}
                \State \textbf{break}
            \EndIf
        \EndWhile    
     \EndFor
    \State \Return $f_{solv}^{down}~\forall solv$   
\end{algorithmic}
\end{algorithm}
\begin{algorithm}[H]
\caption{Feasible lower range}
\begin{algorithmic}[1]
    \For{solv \textbf{in} Pareto}
        \State mode $\gets$ on
        \State counter $\gets$ 0
        \State $cn^{max} \gets cn^{input}$
        \State $cg \gets c_{g,input}$
        
        \While{$cg < cn^{max}$}
            \State $cn \gets cg$ 
            \While{$cn < cn^{max}$}
                \State \textbf{compute} $cost$
                \If{\( [cost_{n}^{wEC} < cost_{n}^{woEC}~\forall n]~\land~[cost_c < 0] \)}
                    \If{mode = off}
                        \State mode $\gets$ on
                        \State \textbf{save} $f_{solv}^{down}(cg_{counter}\gets cg)$
                        \State \textbf{save} $f_{solv}^{down}(cn_{counter,mode} \gets cn)$
    
                        \If{counter = 0}
                            \State counter $\gets$ 1
                        \EndIf
                    \EndIf
                    \State $cn \gets cn + 1$
                \Else
                    \If{mode = on}
                        \State mode $\gets$ off
                        \State \textbf{save} $f_{solv}^{down}(cn_{counter,mode} \gets cn-1)$
                        \State \textbf{break}
                    \EndIf
                    \State $cn \gets cn + 1$
                \EndIf
            \EndWhile
            \State mode $\gets$ off
            \State $cg \gets cg+1$
            \If{counter $> 0$}
                \State counter $\gets$ counter $+1$
            \EndIf
            \If{counter = 6}
                \State \textbf{break}
            \EndIf
        \EndWhile    
    \EndFor
    \State \Return $f_{solv}^{down}~\forall solv$
\end{algorithmic}
\end{algorithm}
\begin{algorithm}[H]
\caption{Find maximum difference between cn and cg}
\begin{algorithmic}[1]
    \For{solv \textbf{in} Pareto}
        \State \textbf{Data:} $f^{up}$
        \State counter $\gets$ 0
        
        \While{counter $< 6$}
            \State $diff_{counter} = cn_{counter}-cg_{counter,mode=off}$ 
            \State counter $\gets$ counter $+1$
        \EndWhile
        \State \textbf{find} $\max(\text{diff})$
        \For{$min[cn_{counter}$ \textbf{with} $diff_{counter} =\max(\text{diff})]$}
            \State \textbf{save} $rang_{solv}e^{up} ($
            \State $spann^{min}\gets cn_{counter}-cg_{counter,mode=on}$
            \State $spann^{max}\gets diff_{counter}$
            \State $cn \gets cn_{counter}$
            \State $cg^{min} \gets cg_{counter,mode=off}$
            \State $cg^{max} \gets cg_{counter,mode=on}$
            \State $)$
        \EndFor
        \State \textbf{Data:} $f^{down}$
        \State counter $\gets$ 0
        \While{counter $< 6$}
            \State $diff_{counter} = cn_{counter,mode=off}-cg_{counter}$ 
            \State counter $\gets$ counter $+1$
        \EndWhile
        \State \textbf{find} $\max(\text{diff})$
        \For{$min[cg_{counter}$ \textbf{with} $diff_{counter} =\max(\text{diff})]$}
            \If{$cg_{counter}$ \textbf{is not in} $range_{solv}^{up}$}
                \State \textbf{save} $range_{solv}^{down} ($
                \State $spann^{min}\gets cn_{counter,mode=on}-cg_{counter}$
                \State $spann^{max}\gets diff_{counter}$
                \State $cg \gets cg_{counter}$
                \State $cn^{min} \gets cg_{counter,mode=on}$
                \State $cn^{max} \gets cg_{counter,mode=off}$
                \State $)$
            \EndIf
        \EndFor
    \EndFor
    \State \Return $range_{solv}^{up},range_{solv}^{down}~\forall solv$
\end{algorithmic}
\end{algorithm}
\begin{algorithm}[H]
\caption{Find minimum difference between cn and cg}
\begin{algorithmic}[1]
    \For{solv \textbf{in} Pareto}
        \State \textbf{Data:} $f^{up}$
        \State counter $\gets$ 0
        
        \While{counter $< 6$}
            \State $diff_{counter} = cn_{counter}-cg_{counter,mode=on}$ 
            \State counter $\gets$ counter $+1$
        \EndWhile
        \State \textbf{find} $\max(\text{diff})$
        \For{$min[cn_{counter}$ \textbf{with} $diff_{counter} =\min(\text{diff})]$}
            \State \textbf{save} $rang_{solv}e^{up} ($
            \State $spann^{min}\gets cn_{counter}-cg_{counter,mode=on}$
            \State $spann^{max}\gets diff_{counter}$
            \State $cn \gets cn_{counter}$
            \State $cg^{min} \gets cg_{counter,mode=off}$
            \State $cg^{max} \gets cg_{counter,mode=on}$
            \State $)$
        \EndFor
        \State \textbf{Data:} $f^{down}$
        \State counter $\gets$ 0
        \While{counter $< 6$}
            \State $diff_{counter} = cn_{counter,mode=on}-cg_{counter}$ 
            \State counter $\gets$ counter $+1$
        \EndWhile
        \State \textbf{find} $\max(\text{diff})$
        \For{$min[cg_{counter}$ \textbf{with} $diff_{counter} =\min(\text{diff})]$}
            \If{$cg_{counter}$ \textbf{is not in} $range_{solv}^{up}$}
                \State \textbf{save} $range_{solv}^{down} ($
                \State $spann^{min}\gets cn_{counter,mode=on}-cg_{counter}$
                \State $spann^{max}\gets diff_{counter}$
                \State $cg \gets cg_{counter}$
                \State $cn^{min} \gets cg_{counter,mode=on}$
                \State $cn^{max} \gets cg_{counter,mode=off}$
                \State $)$
            \EndIf
        \EndFor
    \EndFor
    \State \Return $range_{solv}^{up},range_{solv}^{down}~\forall solv$
\end{algorithmic}
\end{algorithm}
\end{document}